\input gtmacros
\input amsnames
\input amstex
%
\catcode`\@=12        
%
%
\input gtoutput
\volumenumber{3}\papernumber{4}\volumeyear{1999}
\pagenumbers{103}{118}\published{28 May 1999}
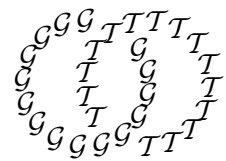
\proposed{Ralph Cohen}\seconded{Haynes Miller, Gunnar Carlsson}
\received{10 May 1998}\revised{5 May 1999}
\accepted{13 May 1999}
%
\let\\\par
\def\topmatter{\relax}
\def\endtopmatter{\maketitlepage}
\let\gttitle\title
\def\title#1\endtitle{\gttitle{#1}}
\let\gtauthor\author
\def\author#1\endauthor{\gtauthor{#1}}
\let\gtaddress\address
\def\address#1\endaddress{\gtaddress{#1}}
\let\gtemail\email
\def\email#1\endemail{\gtemail{#1}}
\def\subjclass#1\endsubjclass{\primaryclass{#1}}
\let\gtkeywords\keywords
\def\keywords#1\endkeywords{\gtkeywords{#1}}
\def\heading#1\endheading{{\def\S##1{\relax}\def\\{\relax\ignorespaces}
    \section{#1}}}
\def\head#1\endhead{\heading#1\endheading}

\def\subhead#1\endsubhead{\sh{#1}}
\def\subsubhead#1\endsubsubhead{\sh{#1}}
\def\specialhead#1\endspecialhead{\sh{#1}}
\def\demo#1{\rk{#1}\ignorespaces}
\def\enddemo{\ppar}
\let\remark\demo
\def\endremark{}
\let\definition\demo
\def\enddefinition{\ppar}

\def\qed{\ifmmode\quad\sq\else\hbox{}\hfill$\sq$\par\goodbreak\rm\fi}  
\def\proclaim#1{\rk{#1}\sl\ignorespaces}
\def\endproclaim{\rm\ppar}
\def\cite#1{[#1]}
\newcount\itemnumber
\def\roster{\items\itemnumber=1}
\def\endroster{\enditems}

\let\itemold\item
\def\item{\itemold{{\rm(\number\itemnumber)}}%
\global\advance\itemnumber by 1\ignorespaces}
\def\S{section~\ignorespaces}  
\def\date#1\enddate{\relax}
\def\thanks#1\endthanks{\relax}   
\def\dedicatory#1\enddedicatory{\relax}  
\let\footnote\plainfootnote

\def\Refs{\ppar{\large\bf References}\ppar\bgroup\leftskip=25pt
\frenchspacing\parskip=3pt plus2pt\small}       
\def\endRefs{\egroup}
\def\widestnumber#1#2{\relax}
\def\endrefitem{}
\def\refdef#1#2#3{\def#1{\leavevmode\unskip\endrefitem#2\def\endrefitem{#3}}} 
\def\ref{\par}
\def\endref{\endrefitem\par\def\endrefitem{}}
\refdef\key{\noindent\llap\bgroup[}{]\ \ \egroup}
\refdef\no{\noindent\llap\bgroup[}{]\ \ \egroup}
\refdef\by{\bf}{\rm, }
\refdef\manyby{\bf}{\rm, }
\refdef\paper{\it}{\rm, }
\refdef\book{\it}{\rm, }
\refdef\jour{}{ }
\refdef\vol{}{ }
\refdef\yr{$(}{)$ }
\refdef\ed{(}{ Editor) }
\refdef\publ{}{ }
\refdef\inbook{from: ``}{'', }
\refdef\pages{}{ }
\refdef\page{}{ }
\refdef\paperinfo{}{ }
\refdef\bookinfo{}{ }
\refdef\publaddr{}{ }
\refdef\eds{(}{ Editors)}
\refdef\bysame{\hbox to 3 em{\hrulefill}\thinspace,}{ }
\refdef\toappear{(to appear)}{ }
\refdef\issue{no.\ }{ }

\define\la{\longrightarrow}

\define\shtimes{\!\times\!}
\define\minus{\smallsetminus}

\define\itb{\itemold{$\bullet$}}

\define\cO{\Cal O}

\define\cX{\Cal X}
\define\cY{\Cal Y}

\define\RR{\Bbb R}
\define\DD{\Bbb D}
\define\SS{\Bbb S}

\define\ZZ{\Bbb Z}

\define\Or{\operatorname{O}}

\define\holim{\operatornamewithlimits{holim}}

\define\hofiber{\operatorname{hofiber}}
\define\imm{\operatorname{imm}}
\define\emb{\operatorname{emb}}
\define\map{\operatorname{map}}
\define\ivmap{\operatorname{ivmap}}
\define\intr{\operatorname{int}}
\topmatter
\title Embeddings from the point of view of\\immersion theory : Part II 
\endtitle
\shorttitle{Embeddings from immersion theory : II}
\author
Thomas G Goodwillie\\Michael Weiss
\endauthor
\address  Department of Mathematics, Brown University\\Box 1917, Providence RI 02912--0001,
USA\endaddress
\secondaddress{Department of Mathematics, University of Aberdeen\\
Aberdeen AB24 3UE, UK}
\asciiaddress{Department of Mathematics, 
Brown University, Box 1917, Providence RI 02912--0001,
USA 
and
Department of Mathematics, University of Aberdeen, 
Aberdeen AB24 3UE, UK}

\email tomg\@math.brown.edu\\m.weiss\@maths.abdn.ac.uk \endemail
\primaryclass{57R40}\secondaryclass{57R42}
\keywords Embedding, immersion, calculus of functors
\endkeywords

\abstract Let $M$ and $N$ be smooth manifolds. For an
open $V\subset M$ let $\emb(V,N)$ be the space of embeddings from
$V$ to $N$.  By the results of Goodwillie 
\cite{4}, \cite{5}, \cite{6} and Goodwillie--Klein \cite{7},
the cofunctor $V\mapsto\emb(V,N)$ 
is {\it analytic} if $\dim(N)-\dim(M)\ge 3$. We deduce that its Taylor 
series converges to it.
For details about the Taylor series, see Part I (reference \cite{15}).
\endabstract

\asciiabstract{Let M and N be smooth manifolds. For an
open V of M let emb(V,N) be the space of embeddings from
V to N.  By results of Goodwillie 
and Goodwillie--Klein, the cofunctor V |--> emb(V,N) 
is analytic if dim(N)-dim(M) > 2. We deduce that its Taylor 
series converges to it. For details about the Taylor series, see Part I.}
\endtopmatter

\newpage

\document
\catcode`\@=\active
\sectionnumber=-1
\section{Introduction}
This is a continuation of \cite{15}. The ideas in this second part are 
mostly due to Goodwillie; notes and illustrations mostly by 
Weiss. For introductions and notation generally
speaking, see \cite{14} and \cite{15}. We fix smooth
manifolds $M^m$ and $N^n$, without boundary unless otherwise stated,
and we assume $m\le n$. As in \cite{15} we  
write $\cO$ for the poset of open subsets of $M$. Then 
$V\mapsto\emb(V,N)$ is a cofunctor from $\cO$ to Spaces (=
fibrant simplicial
sets), and from \cite{14} and \cite{15} we have a fairly thorough 
understanding of its Taylor approximations, the cofunctors
$V\mapsto T_k\emb(V,N)$. Here we show that the Taylor approximations
are good approximations. The main result is Corollary 2.5. Sections 
4 and 5 contain illustrations. 

Following are some conventions. 

Let $R$ be a set with $r$ elements.
A functor $\Cal X$ 
from the poset $\Cal P_R$ of subsets of $R$ to
Spaces is an $r$--dimensional {\sl cube} of Spaces \cite{3, \S1}. 
The cube is 
$k$--{\sl Cartesian} if the canonical map
$$\Cal X(\emptyset)\la\holim_{\{S\mid \emptyset\ne S\subset R\}}\Cal X(S)$$
is $k$--connected. A {\sl Cartesian} cube is one which is $k$--Cartesian for 
all $k$. Since $\Cal P_R$ as a category is isomorphic to its own opposite,
we can use similar terminology for cofunctors
from $\Cal P_R$ to Spaces. Such a cofunctor will also be called a cube.
It is $k$--{\sl Cartesian} if the canonical map
$$\Cal X(R)\la\holim_{\{S\mid S\subset R, S\ne R\}}\Cal X(S)$$
is $k$--connected, and {\sl Cartesian} if it is $k$--Cartesian for all $k$.

Let $\RR^{n+2}_{\,\llcorner}:=[0,\infty)\times[0,\infty)\times\RR^n$.
An $(n+2)$--dimensional {\sl smooth manifold triad} is a 
paracompact Hausdorff space $Q$ 
together with a maximal atlas consisting of open subsets $V_i\subset Q$
and open embeddings $\phi_i\co V_i\to \RR^{n+2}_{\,\llcorner}$ 
satisfying the following conditions:
\roster
\itb The union of the $V_i$ is $Q$. 
\itb The changes of charts $\phi_i(\phi_j)^{-1}$ are smooth where 
defined, and take points in $\RR^{n+2}_{\,\llcorner}$ with vanishing 
first coordinate 
(vanishing second coordinate) to points with vanishing first 
coordinate (vanishing second coordinate).  
\endroster
Let $\partial_0Q$ be the set of those $x\in Q$ which under some 
chart $\phi_i$ map to a point in $\RR^{n+2}_{\,\llcorner}$ with vanishing 
first coordinate. Also, let $\partial_1Q$ be the set of those $x\in Q$ which 
under some chart $\phi_i$ map to a point in $\RR^{n+2}_{\,\llcorner}$ with vanishing 
second coordinate. Then $Q$ is an $(n+2)$--dimensional smooth manifold
with {\sl corners}; its boundary is $\partial_0Q\cup\partial_1Q$ and its 
corner set is $\partial_0Q\cap\partial_1Q$. 

The {\sl handle index} of the manifold triad $Q$ is the smallest integer $a$ such that 
$Q$ can be built from a collar on $\partial_0Q$ by attaching 
handles of index $\le a$. (It may happen that $Q$ {\sl is} a collar on 
$\partial_0Q$.  When that is the case we say that the handle index is $-\infty$.) 

\remark{{Example}} Suppose that $P$ is smooth, with 
boundary, and let $f\co P\to\RR$ be smooth. If $0$ is a regular value
both for $f$ and for $f|\partial P$, then the inverse image of $[0,\infty)$ 
under $f$ is a manifold triad $Q$, with 
$\partial_0Q=\partial P\cap Q$. (If $\partial P=\emptyset$, then of 
course $\partial_0Q=\emptyset$.) 
Any $Q\subset P$ which can be 
obtained in this way will be called a {\sl codimension zero 
subobject} of $P$.
\endremark
 
\section{Excision Estimates} 

Let $Y$ be a smooth $n$--manifold with boundary. Let $Q_0,Q_1,\dots,Q_r$ be smooth compact $n$--manifold 
triads with handle index $q_i$ where $n-q_i\ge3$.
Suppose that smooth embeddings $e_i\co \partial_0Q_i\to\partial Y$
have been selected for $0\le i\le r$, and that their images are 
pairwise disjoint. For $S\subset[r]=\{0,\dots,r\}$ let $Q_S$ be the 
disjoint union of the $Q_i$ with $i\in S$. By $\emb(Q_S,Y)$ we mean 
the Space of smooth embeddings $f\co Q_S\to Y$ which satisfy 
$f\mid\partial_0 Q_i=e_i$ and $f^{-1}(\partial Y)=\partial_0Q_S$. 

\proclaim{1.1 Theorem} 
{\rm(\cite{4}, \cite{5}, \cite{6}, \cite{7})}\qua 
{\sl The $(r+1)$--cube taking a subset $S$ of $[r]$ to 
$\emb(Q_S,Y)$ 
is $(3-n+\Sigma_{i=0}^r(n\!-\!q_i\!-\!2))$--Cartesian, 
provided $r\ge1$.}
\endproclaim

\remark{{Comment}} Assuming that $\emb(Q_{[r]},Y)$ is 
nonempty, fix a base point $f$ in $\emb(Q_{[r]},Y)$. 
Let the image of $f$ in $\emb(Q_S,Y)$ serve as base point 
for $\emb(Q_S,Y)$.  Goodwillie shows in \cite{6} that the cube 
$S\mapsto \Omega\emb(Q_S,Y)$ is 
$(2-n+\Sigma_{i=0}^r(n\!-\!q_i\!-\!2))$--Cartesian, 
provided $r\ge1$.  The delooped statement, 
Theorem 1.1 just above, has been proved by Goodwillie 
and Klein and will appear in \cite{7}. 
\endremark

We will need a slight generalization of 1.1 where the $Q_i$ are 
allowed to have dimension $m\le n$. For this and other purposes 
we need a lemma.

\proclaim{1.2 Lemma} {\sl Let $u\co \cX\to\cY$ be a map of 
$(r+1)$--cubes. That is, $\cX$ and $\cY$ are functors 
from the poset of subsets of $[r]$ to Spaces, and $u$ is a natural 
transformation. Suppose that $\cY$ is $k$--Cartesian and, for 
every $y\in\cY(\emptyset)$, the $(r+1)$--cube defined by 
$$S\mapsto\hofiber[\cX(S)@>u>>\cY(S)]$$
is $k$--Cartesian. {\rm(}The homotopy fiber is taken over the image of 
$y$ in $\cY(S)$.{\rm)} Then $\cX$ is $k$--Cartesian.}
\endproclaim

\demo{Proof of 1.2} Combine \cite{3, 1.18} with 
\cite{3, 1.6}. \qed
\enddemo 

Now let $Y$ and $Q_i$ be as in 1.1, assuming however 
$\dim(Q_i)=m\le n$. As before, embeddings 
$e_i\co \partial_0Q_i\to\partial Y$ are specified and $n-q_i\ge3$, 
where $q_i$ is the handle index of $Q_i$. We want to show that 
the conclusion of 1.1 still holds. Without loss of 
generality, $Y$ is a smooth submanifold of some $\RR^t$. 
Then we can define maps 
$$\emb(Q_S,Y)\la \map(Q_S,G_{n-m})\quad;\quad f\mapsto\nu_f$$
where $G_{n-m}$ is the direct limit of the Grassmannians 
of $(n-m)$--dimensional linear subspaces of $\RR^u$, for $u\ge0$, and 
where $\nu_f$ takes $x\in Q_S$ to the intersection of the tangent space of 
$Y$ at $f(x)$ with  
the normal space of $f(Q_S)$ at $f(x)$. In other words, $\nu_f$ 
classifies the normal bundle of the embedding $f$. We have
therefore a map of $(r+1)$--cubes: 
$$\{S\mapsto \emb(Q_S,Y)\,\}\,\,\la\,\,\{S\mapsto\map(Q_S,G_{n-m})\,\}\,.
\tag$\bullet$ $$
Since the codomain cube in \thetag{$\bullet$} is Cartesian, the contravariant 
analog of 1.2 tells us that the 
domain cube is $k$--Cartesian provided 
that, for every $\xi$ in $\map(Q_{[r]},G_{n-m})$, 
the cube of homotopy fibers 
$$S\mapsto \hofiber_{\xi}\,[\,\emb(Q_S,Y)\}\,\,\la\,\,\map(Q_S,G_{n-m})\,]
\tag{$\bullet\bullet$}$$
is $k$--Cartesian. Now it is easy to construct a homotopy equivalence from
the homotopy fiber of $\emb(Q_S,Y)\to\map(Q_S,G_{n-m})$ to 
$\emb(Q'_S,Y)$ where $Q'_S$ is the total space of the disk bundle 
on $Q_S$ determined by $\xi$. This construction is natural in $S$, so 
\thetag{$\bullet\bullet$} is $k$--Cartesian if the cube
$S\mapsto\emb(Q'_S,Y)$ is $k$--Cartesian, which by 1.1 is the case if
$k=(3-n+\Sigma_{i=0}^r(n\!-\!q_i\!-\!2))$. Therefore: 

\proclaim{1.3 Observation} {\sl Theorem 1.1 generalizes to the situation
where the $Q_i$ have dimension $m\le n$.}
\endproclaim

The corollary below is a technical statement about the
cofunctor on $\cO$ given by $V\mapsto\emb(V,N)$. 
Suppose that $P$ is a smooth compact codimension 
zero subobject of $M$, and that $Q_0,\dots,Q_r$ are
pairwise disjoint compact codimension zero subobjects of 
$M\minus\intr(P)$. For $S\subset\{0,\dots,r\}=[r]$ let 
$V_S$ be the interior of $P\cup Q_S$ where $Q_S$ 
is the union of the $Q_i$ for $i\in S$. We write $V$ instead of $V_{\emptyset}$. 
Suppose that $Q_i$ has handle index $q_i\le n-3$.
 
\proclaim{1.4 Corollary} {\sl Assume $r\ge1$. The
$(r+1)$--cube taking a subset $S$ of $\,[r]\,$ 
to $\,\emb(V_S,N)\,$ is
$(3-n+\Sigma_{i=0}^r(n\!-\!q_i\!-\!2))$--Cartesian.}
\endproclaim

\demo{{\bf Proof}} Let $\bar V_S$ be the closure of $V_S$ in $M$.
Let $\emb(\bar V_S,N)$ be the inverse limit of the Spaces 
(simplicial sets) $\emb(U,N)$ where $U$ ranges over the 
neighborhoods of $V_S$ in $M$. The restriction from 
$\emb(\bar V_S,N)$ to $\emb(V_S,N)$ is a homotopy equivalence. 
Using this fact and 1.2,  we see that it is enough to show 
that for every embedding $f\co \bar V\to N$, the $(r+1)$--cube
$$S\mapsto \hofiber\,[\,\emb(\bar V_S,N)@>\text{res}>>\emb(\bar V,N)\,]$$
is $(3-n+\Sigma_{i=0}^r(n\!-\!q_i\!-\!2))$--Cartesian. The homotopy
fibers are to be taken over the point $f$ in $\emb(\bar V,N)$, which 
we fix for the rest of this proof.
By the isotopy extension theorem, the restriction maps from
$\emb(\bar V_S,N)$ to $\emb(\bar V,N)$ are Kan fibrations. Therefore it is 
enough to show that 
$$S\mapsto \text{fiber}\,[\,\emb(\bar V_S,N)@>\text{res}>>\emb(\bar V,N)\,]$$
is a $(3-n+\Sigma_{i=0}^r(n\!-\!q_i\!-\!2))$--Cartesian $(r+1)$--cube. Let $D(\bar V)$
be the total space of a normal disk bundle for $\bar V$ in $N$, 
with corners rounded off, so that $D(\bar V)$ is a smooth 
codimension zero subobject of $N$. Let $Y$ be the closure of 
$N\minus D(\bar V)$ in $N$. Let $\cX(S)\subset\emb(Q_S,Y)$ 
be the Space of embeddings $g\co Q_S\to Y$ for which the map $f\cup g$ from
$\bar V_S=\bar V\cup Q_S$ to $N$ is smooth. 
(Here all embeddings from $Q_S$ to $Y$ are prescribed on 
$\partial_0Q_S$, as in 1.1 and 1.3.) The inclusions of $\cX(S)$ in 
$\emb(Q_S,Y)$ and the fiber of the restriction map $\emb(\bar V_S,N)@>>>\emb(\bar V,N)$ are 
homotopy equivalences. Therefore it is enough to know that 
$S\mapsto \emb(Q_S,Y)$ is $(3-n+\Sigma_{i=0}^r(n\!-\!q_i\!-\!2))$--Cartesian,
which we know from 1.3. \qed
\enddemo

\section{Convergence}
We begin with an abstraction. 
Suppose that $G$ is a good cofunctor from $\cO$ to Spaces
\cite{14, 2.2}.  Fix an integer $\rho>0$. 
Let $P$ be a smooth compact codimension 
zero subobject of $M$, and let $Q_0,\dots,Q_r$ be
pairwise disjoint compact codimension zero subobjects of 
$M\minus\intr(P)$.  Suppose that $Q_i$ has handle index 
$q_i<\rho$. Let $V_S=\intr(P\cup Q_S)$ as in 1.4. 
 Assume also $r\ge1$.

\definition{2.1 Definition} The cofunctor $G$ is {\sl $\rho$--analytic with 
excess $c$} if,
in these circumstances, the $(r+1)$--cube $S\mapsto G(V_S)$ is
$(c+\Sigma_{i=0}^r(\rho-q_i))$--Cartesian.
\enddefinition

\definition{2.2 Example} According to 1.4, the cofunctor 
$V\mapsto\emb(V,N)$ is $(n-2)$--analytic, with excess $3-n$. 
\enddefinition

\proclaim{2.3 Theorem} {\sl Suppose that $G$ is $\rho$--analytic with excess $c$, 
and that $W\in\cO$ has a proper Morse function whose critical points are all of 
index $\le q$, where $q<\rho$. Then $\eta_{k-1}\co G(W)\to T_{k-1}G(W)$ is 
$(c+k(\rho-q))$--connected, for any $k>1$.} 
\endproclaim

\demo{{\bf Proof}} With a homotopy inverse limit argument we can  
reduce to the case where $W$ is tame (ie, is the interior 
of a compact codimension zero subobject of $M$) and $\bar W$ has a 
smooth handle decomposition with handles of index $\le q$ only.

{\bf Case 1 : $q=0$}\qua Then $W$ is the union of disjoint open  
$m$--balls $W_i$ with $1\le i\le\ell$ for some $\ell$. For 
$S\subset\{1,\dots,\ell\}$ we write $W_S=\cup_{i\in S}W_i$.
If $\ell<k$ then 
$\eta_{k-1}$ from $G(W)$ to $T_{k-1}G(W)$ is a homotopy equivalence
by definition of $T_{k-1}G$. Assume 
therefore $\ell\ge k$. The diagram
$$T_{\ell}G(W)@>r_{\ell}>>T_{\ell-1}G(W)@>r_{\ell-1}>>\cdots
@>r_k>>T_{k-1}G(W)$$
can be identified up to homotopy equivalences with 
$$J_{\ell}@>>>J_{\ell-1}@>>>\cdots@>>>J_k@>>>J_{k-1}$$
where $J_t$ is the homotopy inverse limit of the $G(W_S)$ with 
$S\subset\{1,\dots,\ell\}$ and $|S|\le t$. (See also
the last lines of the proof of \cite{15, 9.1}.) The fibration
$J_t@>>>J_{t-1}$ is obtained by pullback from another fibration, namely 
the product over all $S\subset\{1,\dots,\ell\}$ with $|S|=t$ of the 
fibration
$$p_S\co \holim_{\{R\mid R\subset S\}} G(W_R) \,\,\la\,\,
\holim_{\{R\mid R\subset S, R\ne S\}} G(W_R)\,.$$
All this is true without any special assumptions on $G$ except goodness.
But now we use the analyticity hypothesis 
and find that $p_S$ is $(c+t\rho)$--connected where $t=|S|$. Therefore 
the composition 
$$r_kr_{k+1}\dots r_{\ell}\co T_{\ell}G(W)\to T_{k-1}G(W)$$
is $(c+k\rho)$--connected. We can identify it with 
$\eta_{k-1}\co G(W)\to T_{k-1}G(W)$ since $\eta_{\ell}$ from $G(W)$
to $T_{\ell}G(W)$ is a homotopy equivalence.

{\bf Case 2 : $q>0$}\qua We induct on $q$. For every $q$--handle $Q_u$ 
in the handle decomposition of $\bar W$ choose a smooth chart
$\phi_u\co Q_u\cong\DD^q\times\DD^{m-q}$ and distinct points $x_{u,i}$ in the 
interior of $\DD^q$
for $1\le i\le k$. Let $A_{u,i}\subset W$ be the inverse image of 
$\{x_{u,i}\}\times\intr(\DD^{m-q})$ under $\phi_u$. Each $A_{u,i}$ is 
a closed smooth codimension $q$ submanifold of $W$, meeting the core 
of handle $Q_u$ transversely in one point. For $S\subset\{1,\dots,k\}$
let $A_S$ be the union of all $A_{u,i}$ with $i\in S$ and 
$u$ arbitrary. Our 
analyticity assumption on $G$
implies that the cube defined by $S\mapsto G(W\minus A_S)$ is $(c+k(\rho-q))$--Cartesian.
The cube $S\mapsto T_{k-1}G(W\minus A_S)$ is Cartesian (=$\infty$--Cartesian).
For each nonempty $S$, the manifold $W\minus A_S$ has an isotopy equivalent
tame codimension zero submanifold with a handle decomposition where all 
handles have indices $<q$. Therefore by inductive assumption, the map
$\eta_{k-1}\co G(W\minus A_S)\to T_{k-1}G(W\minus A_S)$
is $(c+k(\rho-q+1))$--connected, provided $S\ne \emptyset$.
It follows with \cite{3, 1.22}
that $\eta_{k-1}$ induces a map
$$\holim G(W\minus A_S)\,\la\,\holim T_{k-1}G(W\minus A_S)$$
which is $(c+k(\rho-q+1)-k+1)$--connected (both homotopy inverse limits 
are over {\sl nonempty} $S\subset\{1,\dots,k\}$). Combining this with our 
``Cartesian--ness'' estimates for the cubes $S\mapsto G(W\minus A_S)$ and 
$S\mapsto T_{k-1}G(W\minus A_S)$, where $S$ again 
denotes an arbitrary subset of $\{1,\dots,k\}$, 
we can conclude that $G(W)\to T_{k-1}G(W)$
is indeed $(c+k(\rho-q))$--connected. \qed
\enddemo

\proclaim{2.4 Corollary} {\sl Suppose that $G$ is $\rho$--analytic. If $\rho>m$,
then the canonical map $G(W)\la\holim_k T_kG(W)$ is a homotopy equivalence 
for every $W$ in $\cO$. In general, the map $G(W)\to\holim_k T_kG(W)$ 
is a homotopy equivalence if $W$ has a proper Morse function whose 
critical points are all of index $<\rho$. }
\endproclaim

\proclaim{2.5 Corollary} {\sl Let $G(W)=\emb(W,N^n)$ for open $W\subset M^m$.
If $m$ is less than $n-2$, then $G(W)\simeq\holim_kT_kG(W)$ for all $W$.
If $m=n-2$, then $G(W)\simeq\holim_kT_kG(W)$ provided 
$W$ has no compact component. 
In general:  suppose that $W$ has a proper Morse function whose critical 
points are all of index $\le q$, where $q<n-2$. Then 
$\eta_k\co G(W)\to T_kG(W)$ is $(k(n-2-q)-q+1)$--connected for $k\ge1$.
Consequently $G(W)\simeq \holim_kT_kG(W)$.}
\endproclaim

\demo{{\bf Proof}} This follows from 2.4, 2.3 and 2.2
with $\rho=n-2$ and $c=3-n$. \qed
\enddemo

\proclaim{2.6 Corollary} {\sl Let $f\co G_1\to G_2$ be a natural 
transformation between good 
cofunctors on $\cO$.  Suppose that $G_1$ and $G_2$ are both $\rho$--analytic and
$f$ from $G_1(W)$ to $G_2(W)$ is a homotopy equivalence whenever $W$ is a 
tubular neighborhood of a finite set ($W\in\cO j$ for some $j$). Then 
$f\co G_1(W)\to G_2(W)$ is a homotopy equivalence for any $W$ which 
has a proper Morse function with critical points of 
index $<\rho$ only.} 
\endproclaim 

\demo{{\bf Proof}} The hypothesis on $f$ implies that $T_kf\co T_kG_1\to T_kG_2$ is an 
equivalence. \qed
\enddemo

\section{Taylor Approximations of Analytic 
Cofunctors} 
As in the preceding section, $G$ is a good cofunctor from $\cO$ to Spaces.

\proclaim{3.1 Proposition} {\sl Suppose that $G$ is homogeneous of degree $k$
where $k\ge0$. If $G(V)$ is $(c-1+k\rho)$--connected for every $V$ 
in $\cO k$, and $\rho\ge m$, then $G$ is $\rho$--analytic 
with excess $c$.}
\endproclaim

\remark{{Remark}} If $V$ in $\cO k$ has $<k$ components, then 
$G(V)$ is contractible and therefore automatically $(c-1+k\rho)$--connected.
The values $G(V)$ for $V$ in $\cO k$ with exactly $k$ components can 
be regarded as the fibers of the {\sl classifying (quasi)--fibration} for $G$.
See \cite{15, \S8} for the classification of homogeneous cofunctors on $\cO$.
\endremark

\demo{{\bf Proof}} The case $k=0$ is trivial, so we assume $k\ge 1$. 
Choose $r$ and $V_S$ for $S\subset[r]$ as in 2.1. 
Again write $V=V_{\emptyset}$. We must show that the cube
$S\mapsto G(V_S)$ is a $(c+\Sigma^r_{i=0}(\rho-q_i))$--Cartesian
$(r+1)$--cube. 
The case $r\ge k$ is again
trivial, so we assume that $r<k$ and proceed by {\sl downward} induction on 
$r$. Also the cases where $q_i=-\infty$ for some $i$ are trivial, so we assume 
$q_i\ge0$ for all $i$. 

{\bf Case 1 :  $q_i=0, \forall i$}\qua We do an induction on the number of 
handles in a handle decomposition of $\bar V$. If $\bar V$ is empty
and $r<k-1$ then $G(V_S)$ is contractible for all $S$, so there is nothing 
to prove. If $\bar V$ is empty and $r=k-1$, then $G(V_S)$ is contractible except 
possibly when $S=[r]$, in which case it is $(c-1+k\rho)$--connected;
then the cube $S\mapsto G(V_S)$ is $(c+k\rho)$--Cartesian, which means,
$(c+\Sigma(\rho-q_i))$--Cartesian. 

For nonempty (but still tame) $\bar V$, choose a handle decomposition.
If $\bar B$ is the cocore of a handle of index $p$ in $\bar V$,
and $B :=\bar B\cap V\cong\RR^{m-p}$, then the $(r+2)$--cube
$$\{G(V_S)\mid S\subset[r]\} \la \{G(V_S\minus B)\mid S\subset[r]\}$$
is $(c+(\rho-p)+(r+1)\rho)$-Cartesian by the downward induction on $r$.
The $(r+1)$--subcube $\{G(V_S\minus B)\mid S\subset[r]\}$ is
$(c+(r+1)\rho)$--Cartesian by the upward induction on the number of
handles of $\bar V$. Noting that $\rho\ge m\ge p$ and using \cite{3, 1.6}
we conclude that the cube $\{G(V_S)\mid S\subset[r]\}$ is also 
$(c+(r+1)\rho)$--Cartesian. 

{\bf Case 2 : $q_0>0$}\qua Let $\bar A$ and $\bar B$ be two parallel but disjoint 
cocores for the handle $Q_0\subset\bar V_{[r]}$ and let $\bar C$ be a strip between 
$\bar A$ and $\bar B$. Let $A=\bar A\cap V$, $B=\bar B\cap V$, $C=\bar C\cap V$
so that the triad $(C;A,B)$ is homeomorphic to the triad
$([0,1]\times\RR^{m-q_0}\,;\,\{0\}\times\RR^{m-q_0}\,,\,\{1\}\times\RR^{m-q_0})$.
Consider the diagram of $r$--cubes

$$\CD
\{G(V_S\minus C)\}@<f<<\{G(V_S\minus(A\cup B))\}@<g<< \{G(V_S\minus B)\} \\
@. @AAA @AAA \\
 @. \{G(V_S\minus A)\}@<h<<\{G(V_S)\}
\endCD$$

where now $S$ runs through subsets of $[r]$ containing the element $0$.
(Each arrow in the diagram is a natural transformation of $r$--cubes,
induced by appropriate inclusions.) 
We have to show that the arrow $h$ (which is a map of $r$--cubes) is a 
$(c+\Sigma(\rho-q_i))$--Cartesian $(r+1)$--cube. It is enough 
\cite{3, 1.6} to show that
$g$ is a $(c+\Sigma(\rho-q_i))$--Cartesian $(r+1)$--cube and the 
square with horizontal arrows $g$ and $h$ is a 
$(c+\Sigma(\rho-q_i))$--Cartesian $(r+2)$--cube. But the composition 
$fg$ is a Cartesian cube, because for each $S$ containing 0 the inclusion
$V_S\minus C\to V_S\minus B$ is an isotopy equivalence. And $f$ is a 
$(c+1+\Sigma(\rho-q_i))$--Cartesian $(r+1)$--cube by induction on $q_0$
(see comments below if this raises doubts). Using \cite{3, 1.8} 
now, we see that $g$ is indeed a $(c+\Sigma(\rho-q_i))$--Cartesian 
$(r+1)$--cube. What about the square? By the downward induction on $r$ 
it is a $(c+(\rho-q_0)+\Sigma(\rho-q_i))$--Cartesian $(r+2)$--cube. This is 
good enough since $\rho\ge m\ge q_0$. \qed
\enddemo

\remark{{Comment}} Let $\bar U$ be the complement of a 
nice tubular neighborhood
of $\bar A\cup\bar B$ in $\bar V_{[r]}$, and $U=\intr(\bar U)$ so that $U$ is tame. 
Then $\bar U$ can be obtained from $\bar V$ by attaching
$r+1$ handles whose indices are $q_0-1,q_1,\dots,q_r$. In fact the 
intersection of $\bar U$ with $\bar C$ is the cocore of the first 
handle in the list, with index $q_0-1$. 

The $(r+1)$--cube determined as in 2.1
by the handle decomposition of $\bar U$ relative to $\bar V$ 
can be identified with the cube $f$ above. More precisely, the two cubes
are ``identified'' by a natural transformation from one to the 
other which is a termwise homotopy equivalence.   
\endremark

\proclaim{3.2 Corollary} {\sl If $G$ is $\rho$--analytic with excess $c$,
then so is $T_kG$.}
\endproclaim

\demo{{Proof}} We proceed by induction on $k$. Fix $V$ and $V_S$ for 
$S\subset[r]$ as in 2.1, with $r\ge0$. 
We must prove that the $(r+1)$--cube $S\mapsto T_kG(V_S)$ is
$(c+\Sigma_{i=0}^r(\rho-q_i))$--Cartesian. For $k=0$ this is correct.
 
Suppose $k\ge1$ and let $z\in T_{k-1}G(V_{[r]})$. For
open $U\subset V_{[r]}$ let  
$L_k^zG(U)$ be the homotopy fiber of the forgetful map 
$T_kG(U)\to T_{k-1}G(U)$ over the point obtained from $z$ 
by restriction. The cofunctor
$U\mapsto L_k^zG(U)$ on open subsets $U\subset V_{[r]}$ is homogeneous of 
degree $k$. If $U\subset V_{[r]}$
is a disjoint union of open $m$--balls, then $L_k^zG(U)$ is 
homotopy equivalent to
the total homotopy fiber of the cube $Y\mapsto G(U_Y)$ where $Y$ 
runs over subsets of $\pi_0(U)$; see \cite{15, 9.1}. By the assumption 
on $G$, this implies that $L_k^zG(U)$ is $(c+k\rho-1)$--connected.
We are now in a position to use 3.1, and conclude that $L_k^zG$ 
is $\rho$--analytic with excess $c$, just like $G$. 
Now the induction step is easy: we can make a 
fibration sequence up to homotopy of $(r+1)$--cubes
$$\{L_k^zG(V_S)\}\la \{T_kG(V_S)\}\la \{T_{k-1}G(V_S)\}\,.$$
By induction, the right--hand cube is $(c+\Sigma(\rho-q_i))$--Cartesian;
the left-hand one also is, for arbitrary $z$,
because we just proved it; and therefore the one in the middle is
$(c+\Sigma(\rho-q_i))$--Cartesian, by lemma 1.3. \qed
\enddemo

\proclaim{3.3 Corollary} {\sl Let $G$ be $\rho$--analytic with excess $c$.
Suppose that $V\in\cO$ has a proper Morse function with critical points of 
index $\le q$ only, where $q<\rho$.
Then $r_k$ from $T_kG(V)$ to $T_{k-1}G(V)$ is 
$(c+k(\rho-q))$--connected, for $k\ge0$.}
\endproclaim

\remark{{Remark}} Note that 3.3 is suggested but not implied by 2.3, since 
$\eta_{k-1}=r_k\eta_k$. Also 3.3 does not imply 2.3 since it only 
depends on the behavior of $G$ on $\bigcup\cO k$. 
\endremark

\demo{Proof of 3.3} We identify $r_k\co T_kG\to T_{k-1}G$ with
$\eta_{k-1}\co T_kG\mapsto T_{k-1}(T_kG)$ (using \cite{15, 6.1}). Now 2.3
can be applied, with $T_kG$ in place of $G$,
because $T_kG$ is $\rho$--analytic with excess $c$ 
according to 3.2. \qed
\enddemo

\section{Haefliger's Theory}
Here we explain how the calculus of embeddings contains the 
Haefliger theory of smooth embeddings in  
the metastable range (which Haefliger unfortunately calls the 
{\sl stable range}). First we recall Haefliger's theory. 
With $M$ and $N$ as before, we can make a commutative diagram

$$\CD
\emb(M,N)@>\subset>>\map(M,N) \\
@VVV   @VV f\mapsto f\times f V \\
\ivmap^{\ZZ/2}(M\shtimes M,N\shtimes N)@>\subset>>
                     \map^{\ZZ/2}(M\shtimes M,N\shtimes N)\endCD\tag{$*$}$$

where all mapping Spaces in sight consist of smooth maps.
Specifically, the expression $\ivmap^{\ZZ/2}(\dots)$ stands for a 
Space of {\sl strictly isovariant} 
smooth maps, ie, equivariant smooth maps $g$ with the properties
$$\align g^{-1}(\Delta_N)=\,&\Delta_M \\
(T_xg)^{-1}(T_{g(x)}\Delta_N)=\,&T_x\Delta_M\quad\text{for any }
x\in\Delta_M\,.\endalign$$ 
(The second of these properties can be reformulated as follows:
the vector bundle morphism induced by $g$, from the normal bundle of 
$\Delta_M$ in $M\shtimes M$ to the normal bundle of $\Delta_N$ in 
$N\shtimes N$, is a fiberwise monomorphism.)

\proclaim{4.1 Theorem}{\rm \cite{8}}\qua {\sl If $2n>3(m+1)$ and $n\ge3$, square 
\thetag{$*$} is
1--Cartesian.}
\endproclaim
 
\remark{{Remark}} In Haefliger's set--up, the space corresponding to 
$\ivmap^{\ZZ/2}(\dots)$ consists of the (smooth  or continuous) 
equivariant maps $g\co M\shtimes M\to N\shtimes N$ for which 
$g^{-1}(\Delta_N)=\Delta_M$. The extra condition on infinitesimal
behavior that we have added does affect the homotopy type. 
However, the comparison map is $1$--connected (proved in
\cite{8} and especially \cite{9}), so that 4.1 is 
correct with either definition. 
\endremark

\proclaim{4.2 Lemma} {\sl If $m\le n$, the cofunctor $E$ taking $V\in\cO$ to the 
homotopy pullback of
$$\ivmap^{\ZZ/2}(V\shtimes V,N\shtimes N)@>\subset>>
            \map^{\ZZ/2}(V\shtimes V,N\shtimes N)@<<<\map(V,N)$$
is polynomial of degree $\le2$.}
\endproclaim

\demo{{\bf Proof}} By \cite{15, 2.5} it suffices to verify that each of the 
three terms from which the homotopy pullback is made is polynomial 
of degree $\le2$, as a cofunctor in $V$. Clearly $V\mapsto\map(V,N)$
is polynomial of degree $\le1$. The remaining two terms can be 
handled as in \cite{15, 2.4}. We omit the details. 
\qed
\enddemo

\proclaim{4.3 Corollary} {\sl Let $F(V)=\emb(V,N)$.
The canonical morphism $F\to E$ is a 
second order Taylor approximation (induces an 
equivalence $T_2F\to T_2E$).}
\endproclaim

\demo{{\bf Proof}} By \cite{15, 5.1} and 4.2 above it suffices to verify that 
$F(V)\to E(V)$
is a homotopy equivalence whenever $V$ is diffeomorphic to a disjoint
union of $\le2$ copies of $\RR^m$. But this is rather obvious. 
\qed
\enddemo

We see that Haefliger's theorem, 4.1, can also be deduced from 4.3 and 2.5.
In fact, 2.5 tells us that $F(M)\to E(M)\simeq  T_2F(M)$ is 
$s$--connected with 
$$s=2(n-2-m)-m+1=2n-3(m+1)\,.$$
Again, this is not a new result. It is a reformulation of the 
main result of Dax's thesis \cite{2}.

\section{An Application}
The convergence statements 2.5 and 4.4 can be generalized mildly so that 
certain cases ``with boundary'' are included. What we have in mind is 
example 10.1 of \cite{15}. Suppose therefore that $M^m$ and $N^n$ 
are smooth manifolds with boundary, $m\le n$, and that a smooth embedding 
$g\co \partial M\to \partial N$ has been selected. 
Let $\cO$ be the poset of open subsets of $M$ 
containing $\partial M$. For $V\in\cO$  let $\emb(V,N)$ be the 
Space of (neat and smooth) embeddings $V\to N$ which agree with $g$
near $V\cap \partial M$. The cofunctor $V\mapsto\emb(V,N)$ on $\cO$ is 
{\sl good}. A calculus of good cofunctors from $\cO$ to 
Spaces is outlined in \S10 of \cite{15}.

\proclaim{5.1 Fact} {\sl Let $F(V)=\emb(V,N)$ for $V\in\cO$.
If $m<n-2$, then $F(V)\simeq\holim_kT_kF(V)$ for all $V$.
In general, $F(V)\simeq\holim_kT_kF(V)$ provided 
$V$ has a handle decomposition relative to a collar on $\partial V$,  with 
possibly infinitely many handles, all of index $\le q$, where $q<n-2$. In this case   
$\eta_k\co F(V)\to T_kF(V)$ is $(k(n-2-q)-q+1)$--connected.}
\endproclaim

The proof is essentially identical with that of 2.5. We omit 
it and turn to some examples. Suppose that $M$ is the unit interval $I$, 
and suppose for simplicity that it comes with a preferred 
smooth and neat embedding $I\to N$ (which we treat as an inclusion).   
Conditions on $N^n$ will be added later; for now, the only condition 
is $n\ge4$. Then $V\mapsto\emb(V,N)$ 
is a functor from $\cO$ to based Spaces. {\sl Notation:}
$\iota X=X\minus\partial X$ for a manifold with boundary $X$.

Define $F(V)$ as in 5.1, for open $V\subset I$ containing the boundary.
Our first and most important task is to understand the 
classifying fibration $p_k$ for $L_kF$, assuming $k\ge2$. Its base space is 
$$\binom{\iota I}k\,\,\cong\,\,\RR^k\,.$$ 
Its fiber over some $S$ (subset of $\iota I$ with $k$ elements) is the total 
homotopy fiber of the $k$--cube of pointed spaces
$$R\mapsto\emb(R,\iota N)\qquad\qquad(R\subset S)\,.\tag{$*$}$$
Let $x\in S$ be the minimal element. In 5.2 and 5.3 below, 
$x\in S\subset I\subset N$ serves as base point of $\iota N$.

\proclaim{5.2 Lemma} {\sl The total homotopy fiber of the $k$--cube \thetag{$*$}
is homotopy equivalent to the total homotopy fiber of the $(k-1)$--cube
of pointed spaces
$$\phantom{\qquad\qquad(R\subset S, x\notin R)}
R\mapsto \iota N\minus R\qquad\qquad
(R\subset S, x\notin R)\,.$$}
\rm

{{\bf Proof}}\qua The idea is to use a Fubini principle:
The total homotopy fiber of the $k$--cube \thetag{$*$} is homotopy equivalent to 
the total homotopy fiber of the $(k-1)$--cube
$$\phantom{xxxxxx}
R\mapsto \hofiber[\emb(R\cup x,\iota N)\to\emb(R,\iota N)]
\qquad\quad(R\subset S, x\notin R)\,.$$
The restriction map $\emb(R\cup x,\iota N)\to\emb(R,\iota N)$ is a fibration, 
so that its homotopy fiber can be replaced by its fiber (over the base 
point), which is $\iota N\minus R$. Note that $R\subset S\subset\iota I
\subset\iota N$.
\qed
\enddemo

\proclaim{5.3 Corollary} {\sl The total homotopy fiber 
of the cube \thetag{$*$} is homotopy equivalent to the
total homotopy fiber of a cube of the form
$$\phantom{xxxxxxxxx}
R\mapsto \iota N\vee(R_+\wedge\SS^{n-1})
\qquad\qquad(R\subset S\minus x)\,.$$}
\rm

{{\bf Remark}}\qua The cube in 5.3 is contravariant like \thetag{$*$}, ie, the 
maps in it are collapsing maps, not inclusion maps. The base point of 
$\iota N$ is still $x$. 
\endremark
 
\demo{{\bf Proof}} Choose a smooth embedding $e\co \iota N\to\iota N$, 
isotopic to the identity, such 
that $e(\iota N)$ and $S\minus x$ are disjoint and $e(x)=x$. For each $y\in S\minus x$
choose a smooth embedding $f_y\co \DD^n\to \iota N$ which maps the base point $(1,0,\dots,0)$
to $x$ and the center $(0,\dots,0)$ to $y$, and avoids all other points of $S$.
Then for each $R\subset S\minus x$, the map
$$\iota N\vee(R_+\wedge\SS^{n-1})\vee((S\minus R\minus x)_+\wedge\DD^n)\la 
\iota N\minus R$$
which is $e$ on $\iota N$ and $f_y$ on $\{y\}_+\wedge\SS^{n-1}$ and on
$\{y\}_+\wedge\DD^n$ (where applicable) is a homotopy equivalence. It is also
contravariantly natural in the variable $R$.
Further, for each $R\subset S\minus x$, the collapse map
$$\iota N\vee(R_+\wedge\SS^{n-1})\vee((S\minus R\minus x)_+\wedge\DD^n)
\la\iota N\vee(R_+\wedge\SS^{n-1})$$
is a homotopy equivalence, again contravariantly natural in the variable $R$.
We have now reduced 5.3 to 5.2. 
\qed
\enddemo

In certain cases the Hilton--Milnor theorem \cite{11} can be used to simplify 
5.3 further. See \cite{16}  
for the meaning of {\sl basic words}, which we use in 5.4 below. 
For a basic word $w$ in the letters $z_1,\dots,z_k$ we let $\alpha(w)$ 
be the number of letters distinct from $z_1$, and $\beta(w)$ the number 
of letters equal to $z_1$.  
For a pointed CW--space $Y$ the smash product  
$Y\wedge Y\dots\wedge Y$ with $j$ factors will be abbreviated $Y^{(j)}$;
the convention for $j=0$ is $Y^{(0)}=\SS^0$.
 
\proclaim{5.4 Corollary} {\sl Suppose that $N\simeq\Sigma Y$, where $Y$ is a 
connected pointed CW--space. Then the loop space of the 
total homotopy fiber of \thetag{$*$} is homotopy equivalent to a
weak product (union of finite products)
$${\prod_w}'\Omega\Sigma^{1+\alpha(w)(n-2)}Y^{(\beta(w))}$$
taken over all basic words $w$ in the letters $z_1,\dots,z_k$ 
which involve all the letters except possibly $z_1$.}
\endproclaim

Before giving the proof, we must restate the Hilton--Milnor theorem.
Let $X_1,\dots,X_k$ be connected 
pointed CW--spaces. For a basic word $w$ in the letters $z_1,\dots,z_k$
denote by $w(X_1,\dots,X_k)$ the space obtained by substituting $X_i$ 
for $z_i$ and $\wedge$ for multiplication; for example, if $k=3$
and $w=(z_2z_1)z_1$, then $w(X_1,X_2,X_3)$ is $(X_2\wedge X_1)\wedge X_1$.

\proclaim{5.5 Hilton--Milnor Theorem}
{\sl $\Omega\Sigma(X_1\vee\dots\vee X_k)\simeq
\prod_w' \Omega\Sigma(w(X_1,\dots,X_k))$, a weak product taken 
over the basic words $w$ in the letters $z_1,\dots,z_k$.}
\endproclaim

\demo{Proof of 5.4} 
We identify $S$ in 5.3 and 5.4 with $\{1,2,\dots,k\}$, noting
that $S$ has a preferred ordering since it is a subset of $I$. Fix $X_1=Y$
as in 5.4, and for $i>1$ regard $X_i$ as a variable with only 
two possible values, $\SS^{n-2}$ and $*$. Letting the $X_i$ vary with 
these constraints, we find that
$$\Omega\Sigma(X_1\vee\dots\vee X_k)$$
runs through the loop spaces of the vertices of the cube in 5.3. Each 
of these loop spaces can therefore be decomposed according to Hilton--Milnor.
Moreover, the decomposition given by Hilton--Milnor is natural, a point
which is stressed in \cite{16}, and we conclude that $\Omega$ of the cube
in 5.3 splits into many 
separate factors indexed by the basic words $w$ in the 
letters $z_1\dots,z_k$. Each factor is a $(k-1)$--cube in its own right,
of the form
$$R\mapsto \Omega\Sigma(w(Y,X_2(R),\dots,X_k(R))\qquad\qquad(R\subset
\{2,\dots,k\})\tag{$**$}$$
where $X_i(R)=\SS^{n-2}$ if $i\in R$ and $X_i(R)=*$ otherwise. If $w$
omits one of the letters $z_2,\dots,z_k$, then the total homotopy 
fiber of \thetag{$**$} is contractible (by the Fubini priciple,
also used in the proof of 5.2). Otherwise,
the cube \thetag{$**$} has only one nontrivial vertex (the initial one)
and its total homotopy fiber is therefore that vertex. \qed
\enddemo

\proclaim{5.6 Corollary} {\sl Suppose again that $N\simeq\Sigma Y$ as in 5.4.
The terms of the Taylor tower of $\emb(I,N)$ are as follows, for $k\ge 2$:
$$L_k\emb(I,N)\simeq{\prod_w}'\Omega^k\Sigma^{1+\alpha(w)(n-2)}Y^{(\beta(w))}$$
where the weak product is over all basic words $w$ in the letters 
$z_1,\dots,z_k$ which involve all the letters except possibly $z_1$.}
\endproclaim

\demo{{\bf Proof}} The $k$--th term of the Taylor tower can be described as the space 
of sections with compact support of a fibration $p_k$. 
See \cite{15, \S10}. We found 
that the base space of $p_k$ is homeomorphic to $\RR^k$, so that the 
section space in question is the $k$--fold loop space of any of 
the fibers. The fibers are described in 5.3, and their loop spaces in 5.4.
Add a prefix $\Omega^{k-1}$ to both sides of the formula in 5.4. \qed
\enddemo

\proclaim{5.7 Summary} {\sl Suppose that $N^n\simeq\Sigma Y$ for a connected pointed 
CW--space $Y$, and $n\ge 4$, and $I$ is neatly embedded in $N$. Then
$\emb(I,N)$ is homotopy equivalent to the homotopy inverse limit of 
a certain diagram
$$\dots\to T_k\emb(I,N)\to T_{k-1}\emb(I,N)\to \dots T_1\emb(I,N)\,.$$
Here $T_1\emb(I,N)$ is $\imm(I,N)$, the Space of smooth immersions which agree 
with the inclusion near $\partial I$. For $k\ge2$, the homotopy fiber of 
the forgetful map $T_k\emb(I,N)\to T_{k-1}\emb(I,N)$ is homotopy equivalent to a weak 
product
$${\prod_w}'\Omega^k\Sigma^{1+\alpha(w)(n-2)}Y^{(\beta(w))}$$
taken over all basic words $w$ in the letters 
$z_1,\dots,z_k$ which involve all the letters except possibly $z_1$.}
\endproclaim

\definition{5.8 Remarks} In 5.7, let $N=\RR^{n-1}\shtimes I$. 
There is a distinguished neat embedding $I\to\RR^{n-1}\shtimes I$ which
identifies $I$ with $0\times I$. There is also a fibration sequence up 
to homotopy 
$$\emb(I,\RR^{n-1}\shtimes I)\la \emb(\SS^1,\SS^n)
@>f>>\Or(n\!+\!1)/\Or(n\!-\!1)\,.$$
The arrow $f$ 
takes a smooth embedding $e\co \SS^1\to\SS^n$ to the orthonormal 2--frame
in $\RR^{n+1}$ consisting of $e(*)$ and the unit tangent vector to 
$e(\SS^1)$ at $e(*)$. Consequently 5.7 calculates $\emb(\SS^1,\SS^n)$ if 
$n\ge 4$, up to extension problems in the homotopy category. We do not know what the 
extensions are, but we do hope that homotopy theorists will be intrigued. 

Unfortunately our convergence result does not cover the approximation map 
from $\emb(I,\RR^2\times I)$ to $\holim_k\,T_k\emb(I,\RR^2\times I)$. 
With or without convergence, we get a map of component sets 
$$\pi_0\emb(I,\RR^2\times I)\to \lim_k\pi_0T_k \emb(I,\RR^2\times I)\,.$$
It is generally believed \cite{1} that this is closely related 
to Vassiliev's theory of knot invariants \cite{12}, \cite{13}. 
Conversely, Kontsevich \cite{10}
used ideas related to Vassiliev's theory to prove results 
similar to 5.7, more specifically, to set up a spectral sequence 
converging to the rational cohomology of 
$\emb(\SS^1,\RR^n)$ for $n\ge4$.  He was able to show 
that the spectral sequence collapses at the $E_2$ term.  One 
hopes that this result can be reformulated and perhaps even 
re--proved in calculus language. 
\enddefinition

\definition{Acknowledgments} It is a pleasure to thank John Klein for 
discussions related to the material in \S5 and for pointing out 
an error in an earlier version of \S2. 

The authors are partially
supported by the NSF.
\enddefinition\eject

\enddocument